\documentclass[12pt]{amsart}
\usepackage{amsmath,amssymb,amsfonts,amsthm,amsopn,dsfont}
\usepackage[dvips]{graphicx}
\usepackage{srcltx}
\usepackage{setspace}
\newtheorem{theorem}{Theorem}[section]
\newtheorem{corollary}[theorem]{Corollary}

\newtheorem{example}[theorem]{Example}

\newtheorem{lemma}[theorem]{Lemma}
\numberwithin{equation}{section}

\newtheorem{proposition}[theorem]{Proposition}
\newtheorem{remark}[theorem]{Remark}

\renewcommand{\ell}{l}
\renewcommand{\epsilon}{\varepsilon}

\def\rd{\bR^d}

\def\R{\right)}

\def\<{\left<}
\def\>{\right>}

\def\mv1{M_v^1}

\def\mn{(m,n)}
\def\mn'{(m',n')}

\newcommand{\usp}{\|u(\sigma)\|^{p-1}_{L^\infty}}
\newcommand{\utp}{\|u(\tau)\|^{p-1}_{L^\infty}}

\def\N{\mathbb{N}}
\def\R{\mathbb{R}}
\def\C{\mathbb{C}}
\def\rd{\mathbb{R}^d}

\begin{document}

\title[On the radius of spatial analyticity for symmetric systems]{On the radius of spatial analyticity for semilinear symmetric hyperbolic systems}
\begin{abstract} 
We study the problem of propagation of analytic regularity for semi-linear symmetric hyperbolic systems. We adopt a global perspective and we prove that if the initial datum extends to a holomorphic function in a strip of radius (=width) $\epsilon_0$, the same happens for the solution $u(t,\cdot)$ for a certain radius $\epsilon(t)$, as long as the solution exists. Our focus is on precise lower bounds on the spatial radius of analyticity $\epsilon(t)$ as $t$ grows. \par
We also get similar results for the Schr\"odinger equation with a real-analytic electromagnetic potential.
\end{abstract}
\author{Marco Cappiello, Piero D'Ancona \and Fabio Nicola}
\address{Dipartimento di Matematica,  Universit\`{a} degli Studi di Torino,
Via Carlo Alberto 10, 10123
Torino, Italy}
\address{Sapienza -- Universit\`a di Roma, Dipartimento di Matematica,
Piazzale A. Moro 2, I-00185 Roma, Italy}
\address{Dipartimento di Scienze Matematiche, Politecnico di
Torino, Corso Duca degli
Abruzzi 24, 10129 Torino,
Italy}
\email{marco.cappiello@unito.it}
\email{dancona@mat.uniroma1.it}
\email{fabio.nicola@polito.it}
\date{}
\maketitle

\section{Introduction}
The main concern in this paper is the propagation of the analytic regularity for systems of semilinear evolution equations. Local propagation of the analyticity for nonlinear partial differential equations has been studied in several papers, see for instance \cite{AM, BG1, BG2,ds1,ds2}, starting from classical results of existence and uniqueness. 
In our paper we focus our attention on systems of semilinear equations with coefficients defined for $(t,x) \in \overline{\R}_+ \times \rd,$ hence globally defined in the space variables. Motivation for this type of study comes from the global analytic theory of nonlinear evolution PDE started with the paper by Kato and Masuda \cite{KM} and based on the proof of the existence of global analytic solutions in the space variables for the related Cauchy problem. After \cite{KM}, many authors proved results of this type for several models as the (generalized) Korteweg-de Vries equation, the Euler equations, the Benjamin-Ono equation, the nonlinear Schr\"odinger equation, see for instance \cite{benjamin, bo1, ghhp, gru, hhg, hay1, hay2, hay3,hp}. A parallel study in an elliptic setting devoted in particular to the analyticity of travelling waves has been developed in \cite{A3, BL1, BL2, CGR2,  CN1}. Besides the existence and uniqueness of an analytic solution, also the estimate of the radius of analyticity represents an interesting issue. Although there exist results in the literature treating particular models, cf. \cite{BK, bo2, KV}, we are not aware of similar results for general semilinear hyperbolic systems on $\rd.$ \\
In the present paper we consider systems with initial data analytic on $\rd$ and admitting a holomorphic extension in a strip of the form $\{x+iy: |y|< \epsilon_0\}$ for some $\epsilon_0 >0.$ We assume moreover the existence of a smooth ($C^{\infty}$) solution of the Cauchy problem. Under these conditions, we prove that the solution is indeed analytic for every $t$ in a strip of the form $\{x+iy \in \C^d: |y|<\epsilon(t)\}$ and we give a precise estimate of lower bounds for $\epsilon (t)$ as $t$ grows. The class of systems to which our results apply includes semilinear symmetric hyperbolic systems but also strictly hyperbolic systems, see Remark \ref{strictly}, and other types of equations like the Schr\"odinger equation with electromagnetic potential (cf.\ \cite{dancona}). For the sake of simplicity, in this Introduction we shall consider the special case of differential operators; we refer to Theorem \ref{mainteo2} below for a more general statement involving pseudodifferential operators. \par
Consider a system of the form 
\[
Lu=\mathcal{N}[u]
\]
in $\overline{\R}_+\times \rd\ni(t,x)$, $u=(u_1,\ldots u_n)$, with 
\begin{equation}\label{equazione}
L=\partial_t+iA_0(i\nabla_x)+\sum_{j=1}^d A_j(t,x)\partial_j +B(t,x),
\end{equation}
where $A_0(i\nabla_x)$ is any formally self-adjoint Fourier multiplier (i.e.\ the symbol $A_0(\xi)$ is a smooth Hermitian $n\times n$ matrix), and $A_j(t,x)$ are $n\times n$ Hermitian matrices with entries in $C(\overline{\R}_+;C^\infty_b(\mathbb{R}^d))$, analytic with respect to $x$, satisfying  
\begin{equation}\label{ipo1}
\sup_{(t,x)\in\overline{\R}_+\times\rd}\|\partial^\alpha_x A_j(t,x)\|\leq C^{|\alpha|+1} \alpha!,\quad \alpha\in\N^d
\end{equation}
for some constant $C>0$ (here $C^\infty_b(\mathbb{R}^d)$ denotes the space of all $C^{\infty}$ functions on $\rd$ bounded with all their derivatives). We assume moreover that $B(t,x)$ is a $n\times n$ matrix with entries in $C(\overline{\R}_+;C^\infty_b(\mathbb{R}^d))$, satisfying 
\begin{equation}\label{1.2bis}
\sup_{(t,x)\in\overline{\R}_+\times\rd}\|\partial^\alpha_x B(t,x)\|\leq C^{|\alpha|+1} \alpha!,\quad\alpha\in\N^d.
\end{equation}
The components of the non-linearity $\mathcal{N}[u]$ are assumed to be polynomials of degree $p\in\N$, $p\geq 2$, in the components of $u$ with analytic coefficients, namely:
\begin{equation}\label{ipo10}
\mathcal{N}[u]_k=\sum_{\gamma\in\N^n:\,2\leq |\gamma|\leq p}g_{k,\gamma}(t,x) u^\gamma,\quad k=1,\ldots,n,
\end{equation}
with $g_{k,\gamma}(t,x)$ continuous, satisfying 
\begin{equation}\label{ipo1bis}
\|\partial^\alpha_x g_{k,\gamma}\|_{L^\infty(\overline{\R}_+\times\R^d)}\leq M^{|\alpha|+1} \alpha!,\quad \alpha\in\N^d
\end{equation}
for some constant $M>0$.\par
When $A_0=0$ in \eqref{equazione} we get what is known as a {\it symmetric hyperbolic system} with bounded coefficients. However, it is important to allow the skew-adjoint term $iA_0(i\nabla_x)$ too, in view of the applications to Schr\"odinger equation below. Incidentally, we observe that the form of this term is admittedly very special, but the reader will note from the proof that allowing more general skew-adjoint operators entails challenging issues.  \par 
We study the above system in the space of uniformly analytic functions, defined as follows.\par
Let $\mathcal{A}$ be the space of functions $f$ satisfying the estimates 
\[
\|\partial^\alpha f\|_{L^\infty}\leq C^{|\alpha|+1}\alpha!,\quad \alpha\in\N^d,
\]
for some constant $C>0$. Every function $f \in \mathcal{A}$ extends to a holomorphic function $f(x+iy)$ in the strip $\{x+iy\in\C^d:\ |y|<C^{-1}\}$ and we say that it has radius of analiticity $\geq C^{-1}$. For technical reasons it is convenient to work with an equivalent definition, where the $L^\infty$ norm is replaced by the norm $\|\cdot\|_s$ in the Sobolev space $H^s(\rd)$ for some $s>d/2$, fixed once and for all.\par
By Sobolev embedding theorem it is easy to prove that $\mathcal{A}=\cup_{\epsilon>0}\mathcal{A}(\epsilon)$, where $\mathcal{A}(\epsilon)$ is given by the smooth functions $f\in L^2(\rd)$ such that $\sup_N E^\epsilon_{N}[f]<\infty$, with 
\[
E^\epsilon_{N}[f]:=\sum_{|\alpha|=N}\frac{\epsilon^{|\alpha|-1}}{\alpha!}\|\partial^\alpha f\|_s,\quad N\geq 1. 
\]  

We have then the following result. 
\begin{theorem}\label{mainteo}
There exists a positive number $\overline{\epsilon}_0>0$, depending only on the equation, such that for all $0<\epsilon_0\leq \overline{\epsilon}_0$ the following result holds true. \par
Let $0\in I\subset[0,+\infty)$ be a time interval. Let $u_0\in\mathcal{A}(\epsilon_0)$, and  assume that there exists $s_0 \geq 0$ such that the Cauchy problem $$\begin{cases} Lu=\mathcal{N}[u] \qquad (t,x) \in I \times \rd, \\ u(0,x) =u_0(x) \end{cases}$$ admits a solution $u\in C^1(I,H^{s'}(\rd)) \cap  C^0(I,H^{{s'}+1}(\rd))$ for every $s' \geq s_0$. Then 
$u(t,\cdot)\in\mathcal{A}(\epsilon(t))$, for every $t\in I$, with
\begin{equation}\label{epst}
\epsilon(t)=\epsilon_0\exp\big(-A\int_0^t (1+\|u(\sigma)\|^{p-1}_{L^\infty})\,d\sigma\big),
\end{equation}
for some constant $A>0$ depending only on $\|u_0\|_{s}+\sup_{N\geq 1}E^{\epsilon_0}_{N}[u_0]$. When the equation is linear $(p=1)$, the constant $A$ depends only on the equation and it is independent of $u_0$. 
\end{theorem}
It will follow from the proof that we can take $A\simeq (\|u_0\|_{s}+\sup_{N\geq 1}E^{\epsilon_0}_{N}[u_0])^{p-1}$.  \par
The threshold $\overline{\epsilon}_0$ corresponds, broadly speaking, to the radius of analyticity of the coefficients $A_j$, $B$ and $g_{k,\gamma}$ in the equation (see \eqref{ipo1}, \eqref{1.2bis} and \eqref{ipo1bis}). Hence $\mathcal{A}(\overline{\epsilon}_0)$ is in some sense the larger space preserved by the equation.\par
Roughly speaking the phenomenon of shrinking of the radius of spatial analyticity is due to the combination of two effects: for linear equations the solution may already compress more and more as $t$ grows; this produces a growth of the spatial derivatives and causes the exponential decay in $\eqref{epst}$. In the presence of a non-linearity this growth may be amplified (if the initial datum is not too small), and this causes the presence of the integral term in \eqref{epst}. We refer to Section \ref{osservazioni} for two model examples which illustrate these facts. Some variants are also discussed there, together with speculations on other lower bounds for $\epsilon(t)$, as well as the issue of the dependence of the constants on the initial datum.\par
As anticipated, the above result applies also to the Schr\"odinger equation with real analytic electromagnetic potential. Namely, we have the following result.\par 
 \begin{corollary}\label{corollario}
Let $A(t,x)=(a_1(t,x),\ldots,a_d(t,x))$ and $V(t,x)$ be the magnetic and the electric potential, respectively, with $a_j(t,x)$ and $V(t,x)$ in $C(\overline{\R}_+;C^\infty_b(\mathbb{R}^d))$, with $a_j(t,x)$ real-valued, and satisfying
\[
\|\partial^\alpha_x a_j\|_{L^\infty(\mathbb{R}_+\times\rd)}\leq C^{|\alpha|+1}\alpha!,\ \alpha\in\N^d,\qquad \|\partial^\alpha_x V\|_{L^\infty(\mathbb{R}_+\times\rd)}\leq C^{|\alpha|+1}\alpha!,\ \alpha\in\N^d.
\]
Consider the magnetic Laplacian 
\[ 
\Delta_A=\sum_{j=1}^d (\partial_j+ia_j)^2=\Delta+2iA\cdot\nabla+i{\rm div}\, A-\sum_{j=1}^d a_j^2
 \]
 and the equation 
 \begin{equation}\label{magnpot}
Lu= \partial_tu-i\Delta_Au-iV(t,x)u=\mathcal{N}[u],
\end{equation}
where the non-linearity $\mathcal{N}[u]$ has the same form as in \eqref{ipo10} and \eqref{ipo1bis} (without the subscript $k$, because we are now considering a scalar case). \par
Then the result of Theorem \ref{mainteo} holds true for the equation \eqref{magnpot}.
\end{corollary}
In fact, it is sufficient to apply Theorem \ref{mainteo} with $A_0(i\nabla_x)=-\Delta$, $A_j(t,x)=2a_j(t,x)$ (which are assumed real-valued, hence Hermitian as $1\times 1$ matrices) and $B={\rm div}\, A+i\sum_{j=1}^d a_j^2+iV$ in \eqref{equazione}.\par
 
\vskip0.3cm 
The core of the proof of Theorem \ref{mainteo} and of the more general Theorem \ref{mainteo2} below consists in proving a suitable analytic energy estimate for the solution of the Cauchy problem. In future works we plan to prove similar results for more general nonlinear terms, possibly treating the case of quasi-linear or fully nonlinear systems, and also to improve the above lower bound for the radius of analyticity for special classes of equations (e.g.\ with constant coefficients). We also recall that in some recent papers \cite{CN1, CN2, CN3}
we proved holomorphic extensions of the solutions of elliptic equations in conical subsets of $\C^d$ of the form $\{x+iy \in \C^d : |y|< \varepsilon (1+|x|)\}, \epsilon >0,$ improving some of the results of \cite{BL1, BL2, CGR2} on the analyticity in a strip for travelling waves and stationary Schr\"odinger equations. An interesting issue would be to investigate the same problem for the corresponding evolution equations and to estimate the decay of the width of the sector as the time variable grows.

 \section{Statement and proof of the main result}
 Let us now consider a more general system of the form
\[
Lu=\mathcal{N}[u]
\]
in $\overline{\R}_+\times \rd\ni(t,x)$, $u=(u_1,\ldots,u_n)$, with 
\begin{equation}\label{equazione-bis}
L=\partial_t+iA_0(D_x)+A(t,x,D_x),
\end{equation}
where $A_0(D_x)$ is any formally self-adjoint Fourier multiplier (i.e.\ the symbol $A_0(\xi)$ is a smooth Hermitian $n\times n$ matrix), and $A(t,x,D_x)$ is a system of pseudodifferential operators, i.e.\
\[
A(t,x,D_x)u=(2\pi)^{-d}\int_{\rd}e^{ix\xi}A(t,x,\xi)\hat{u}(\xi)\,d\xi
\]
whose symbol $A(t,x,\xi)$ is a $n\times n$ matrix satisfying the following conditions: 
\[
\langle\xi\rangle^{|\beta|-1}\partial^\alpha_x\partial^\beta_\xi A(t,x,\xi)\ \textrm{has entries in}\ C(\overline{\R}_+;C^0_b(\mathbb{R}^d\times\R^d))\quad \forall\alpha,\beta\in\N^d
\]
with bounds
\begin{equation}\label{simbolo}
\|\partial^\alpha_x\partial^\beta_\xi A(t,x,\xi)\|\leq C_\beta C^{|\alpha|+1}\alpha!\langle \xi\rangle^{1-|\beta|}, \quad \alpha,\beta\in\N^d\end{equation}
and moreover
\begin{equation}\label{simbolo2}
\sup_{t\in\overline{\R}_+;x,\xi\in\rd}\|A(t,x,\xi)+A(t,x,\xi)^\ast\|<\infty.
\end{equation}
It is immediate to see that the simpler case considered in the Introduction, namely the operator $L$ in \eqref{equazione}, \eqref{ipo1}, \eqref{1.2bis} is in fact a special case. \par
 The non-linearity $\mathcal{N}[u]$ is the same considered in the Introduction, therefore is given by \eqref{ipo10}, with coefficients $g_{k,\gamma}(t,x)$ satisfying \eqref{ipo1bis}. 

\begin{theorem}\label{mainteo2} Let $L$ be given by \eqref{equazione-bis}, \eqref{simbolo}, \eqref{simbolo2} and $\mathcal{N}[u]$ as in \eqref{ipo10}, \eqref{ipo1bis}. Then the same result as in Theorem \ref{mainteo} holds true. 
\end{theorem}

Theorem \ref{mainteo2} will follow if we prove that $E^{\epsilon(t)}_{N}[u(t)]$, $N\geq 1$, is bounded, with $\epsilon(t)$ as in \eqref{epst}. Namely, we will prove by induction that for some constant $C_0>0$ we have
\begin{equation}\label{tesi}
E^{\epsilon(t)}_{N}[u(t)]\leq \Phi(t):=
 C_0 \exp\big(C_0\int_0^t (1+\usp)\,d\sigma \big),
\end{equation}
for every $N\geq 1$, if $\epsilon_0\leq \overline{\epsilon}_0$, with $\overline{\epsilon}_0$ small enough, and the constant $A$ in \eqref{epst} is large enough. 
\par
Under the hypothesis of Theorem \ref{mainteo2} we have the following classical energy estimates 
\begin{equation}\label{energy}
\|v(t)\|_s\leq C e^{Ct} \|v(0)\|_s+ C e^{Ct}\int_0^t e^{-C\sigma}\|Lv(\sigma)\|_s\,d\sigma
\end{equation}
for every function $v\in C^1(I,H^s)\cap C^0(I,H^{s+1})$ and $t\in I$, where $I\subset[0,+\infty)$ is a time interval containing $0$, and the constant $C>0$ depends on $L$ and $s$, see \cite[Proposition 6.1]{Mi}. Usually this estimate is proved without the term $iA_0(i\nabla_x)$ in \eqref{equazione}, but since it does not depend on $x$ and is skew-adjoint, it does not affect the classical proof (in particular the fact that it commutes with spatial derivatives is useful to extend the $L^2$-energy estimates to the $H^s$ version); e.g.\ the proof of \cite[Proposition 2.1.2]{rauch} can be repeated almost verbatim for our operator $L$ in \eqref{equazione-bis}. \par
Now, we apply the above estimate to $v=\partial^\alpha u$, where $u$ solves $Lu=\mathcal{N}[u]$, $u(0)=u_0$. We get 
\[
\|\partial^\alpha u(t)\|_s\leq C e^{Ct} \left(\|\partial^\alpha u_0\|_s+ \int_0^t e^{-C\sigma}(\|[L,\partial^\alpha]u(\sigma)\|_s+\|\partial^{\alpha} \mathcal{N}[u(\sigma)]\|_s)\,d\sigma \right).
\]
We multiply by $\epsilon(t)^{|\alpha|-1}/\alpha!$, $|\alpha|=N\geq1$, and we obtain, since $\epsilon(t)\leq\epsilon(0)=\epsilon_0$, 
\begin{multline}\label{astast}
e^{-Ct}E^{\epsilon(t)}_{N}[u(t)]\leq C E^{\epsilon_0}_{N}[u_0]+C \int_0^t \frac{\epsilon(t)^{N-1}}{\epsilon(\sigma)^{N-1}}e^{-C\sigma}\sum_{|\alpha|=N}\frac{\epsilon(\sigma)^{|\alpha|-1}}{\alpha!}\| [L,\partial^\alpha]u(\sigma)\|_s\,d\sigma\\ +C \int_0^t \frac{\epsilon(t)^{N-1}}{\epsilon(\sigma)^{N-1}}e^{-C\sigma}E^{\epsilon(\sigma)}_N[\mathcal{N}[u(\sigma)]]\,d\sigma.
\end{multline} 
Now, we can choose the constant $C_0$ in \eqref{tesi} such that 
\begin{equation}\label{posizione}
C_0>C \max\{4\sup_NE^{\epsilon_0}_{N}[u_0],2\}
\end{equation}
where $C$ is the constant arising in \eqref{energy}. Observe that by the condition \eqref{posizione} we have $C<C_0$ and $C\sup_NE^{\epsilon_0}_{N}[u_0]<C_0/4$ so that the first term in the right-hand side of \eqref{astast} is then dominated by $\frac{1}{4}e^{-Ct}\Phi(t)$.  We will then prove that \eqref{tesi} holds for every $N\geq 1$ possibly for a bigger constant $C_0$ (independent of $N$). More precisely, we will prove directly the case $N=1$ whereas for $N>1$ we shall argue by induction supposing that \eqref{tesi} holds for subscripts $<N$. Under this hypothesis we will verify in the next sections the following estimates on the two integrals involved in the right-hand side of \eqref{astast}. \par

\begin{proposition}\label{procomm} There exist constants $C'>0$ and $\overline{\epsilon}_0>0$, depending only on the equation, such that, for every $\epsilon_0\leq \overline{\epsilon}_0$  and every $A\geq 1$ in \eqref{epst}, and for every sufficiently large constant $C_0$ in \eqref{tesi} we have
\begin{multline}\label{comm} 
\int_0^t \frac{\epsilon(t)^{N-1}}{\epsilon(\sigma)^{N-1}}e^{-C\sigma}\sum_{|\alpha|=N}\frac{\epsilon(\sigma)^{|\alpha|-1}}{\alpha!}\|[L,\partial^\alpha]u(\sigma)\|_s\,d\sigma\\
 \leq \int_0^t \frac{\epsilon(t)^{N-1}}{\epsilon(\sigma)^{N-1}}C'N e^{-C\sigma}E^{\epsilon(\sigma)}_N[u(\sigma)]\,d\sigma+\frac{1}{4}e^{-Ct}\Phi(t)
\end{multline}
where $\Phi(t)$ is defined in \eqref{tesi}.
\end{proposition}

\begin{proposition}\label{prononlin} There exist constants $C'>0$ and $\overline{\epsilon}_0>0$, depending only on the equation, such that, for every $\epsilon_0\leq \overline{\epsilon}_0$ and for every sufficiently large constants $A$ in \eqref{epst} and $C_0$ in \eqref{tesi}, we have
\begin{multline}\label{nonlin} 
\int_0^t \frac{\epsilon(t)^{N-1}}{\epsilon(\sigma)^{N-1}}e^{-C\sigma}E^{\epsilon(\sigma)}_N[\mathcal{N}[u(\sigma)]]\,d\sigma\\ 
\leq \int_0^t \frac{\epsilon(t)^{N-1}}{\epsilon(\sigma)^{N-1}}C'\usp e^{-C\sigma}E^{\epsilon(\sigma)}_N[u(\sigma)]\,d\sigma+\frac{1}{4}e^{-C t}\Phi(t)
\end{multline}
where $\Phi(t)$ is defined in \eqref{tesi}.
\end{proposition}

Assume that Propositions \ref{procomm} and \ref{prononlin} hold true and set $\Psi_N(t)=e^{-Ct}E^{\epsilon(t)}_{N}[u(t)].$ We can continue the computation in \eqref{astast} as 
\begin{equation}\label{astastast}
\Psi_N(t)\leq \frac{3}{4}e^{-Ct}\Phi(t)+\epsilon(t)^{N-1} \int_0^t \frac{\Psi_N(\sigma)}{\epsilon(\sigma)^{N-1}} C'(N+\usp)\,d\sigma.
\end{equation}
We now use the following form of Gronwall inquality, which can be deduced easily from the classical one (see e.g.\ \cite[Lemma 2.1.3]{rauch}).
\begin{lemma}\label{gronwall}
Let $0\leq g,\psi,a\in L^\infty_{loc}([0,T])$, with $a>0$, and $0\leq h\in L^1_{loc}([0,T])$, and 
\[
\psi(t)\leq g(t)+a(t)\int_0^t h(\sigma)\frac{\psi(\sigma)}{a(\sigma)}\,d\sigma.
\]
Then, with $H(t):=\int_0^t h(\sigma)\, d\sigma$, 
\[
\psi(t)\leq g(t)+a(t) e^{H(t)}\int_0^t e^{-H(\sigma)}h(\sigma)\frac{g(\sigma)}{a(\sigma)}\,d\sigma.
\]
\end{lemma}
We apply Lemma \ref{gronwall} to \eqref{astastast} with $\psi(t)=\Psi_N(t)$, $g(t)=\frac{3}{4}e^{-C t}\Phi(t)$, $h(\sigma)=C'(N+\usp)$, $a(t)=\epsilon(t)^{N-1}$. \par
From \eqref{astastast} we therefore get 
\[
\Psi_N(t)\leq \frac{3}{4}e^{-C t}\Phi(t)+\tilde{\Phi}_N(t)
\]
with 
\begin{multline*}
\tilde{\Phi}_N(t)=\frac{3}{4}\epsilon(t)^{N-1} e^{C'Nt+C'\int_0^t\usp}\int_0^t e^{-C'N\sigma-C'\int_0^\sigma \utp\, d\tau} C'(N+\usp)\\
\times\frac{e^{-C \sigma}\Phi(\sigma)}{\epsilon(\sigma)^{N-1}}\,d\sigma.
\end{multline*}
We want to prove that 
\begin{equation}\label{daver}
\tilde{\Phi}_N(t)\leq \frac{1}{4}e^{-Ct}\Phi(t),
\end{equation} 
if the constants $A$ in \eqref{epst} and $C_0$ in \eqref{tesi} are large enough, so that we get $\Psi_N(t)\leq e^{-C t}\Phi(t)$ and the proof of \eqref{tesi} is concluded.\par To this end we use the explicit expression of $\epsilon(t)$ in \eqref{epst}; we have
\begin{multline*}
\tilde{\Phi}_N(t)=\frac{3}{4}C_0 C'N  e^{(-A(N-1)+C'N)t+(C'-A(N-1))\int_0^t\usp d\sigma}\\
\times
\int_0^t e^{(-C'N+A(N-1)-C+C_0)\sigma+(-C'+C_0+A(N-1))\int_0^\sigma \utp\, d\tau}  d\sigma\\
+ \frac{3}{4}C_0 C' e^{(-A(N-1)+C'N)t+(C'-A(N-1))\int_0^t\usp\,d\sigma}\\
\times \int_0^t e^{(-C'N+A(N-1)-C+C_0)\sigma+(-C'+C_0+A(N-1))\int_0^\sigma \utp\, d\tau}\usp d\sigma.
\end{multline*}
In the first integral we estimate 
\[
e^{(-C'+C_0+A(N-1))\int_0^\sigma \utp\, d\tau}\leq e^{(-C'+C_0+A(N-1))\int_0^t \utp\, d\tau},
\]
which holds if $C_0> C'$, whereas in the second integral we use 
\[
e^{(-C'N+A(N-1)-C+C_0)\sigma}\leq e^{(-C'N+A(N-1)-C+C_0)t},
\] 
which holds if $C_0> C'+C$, $A> 2C'$. We get 
\[
\tilde{\Phi}_N(t)\leq \frac{3}{4}\Big[\frac{C'N}{-C'N+A(N-1)-C+C_0}+\frac{C'}{-C'+C_0+A(N-1)}  \Big] e^{-Ct}\Phi(t)
\]
We conclude by observing that for $N=1$ the expression in parenthesis is $<1/3$ if $C_0$ is large enough, while the same holds for every $N\geq2$, if $A$ is large enough. This shows that \eqref{daver} holds true and concludes the proof. \qed

\section{Proof of Proposition \ref{procomm}}
We have, for $k=1,\ldots,n$,
\[
([L,\partial^\alpha]u)_k=-\sum_{h=1}^n\sum_{0\not=\gamma\leq\alpha}\binom{\alpha}{\gamma}(\partial^\gamma A)(t,x,D_x)_{k,h}\partial^{\alpha-\gamma} u_h,
\]
In view of \eqref{simbolo} the symbols $\partial^\gamma_x A(t,x,\xi)_{k,h}$ belong to the H\"ormander's classes $S^1_{1,0}$ and therefore the corresponding operators are bounded $H^{s+1}\to H^{s}$ with operator norm dominated by a seminorm of their symbol in that symbol class, cf. (\cite[Chapter XVIII]{hormanderIII}). Hence, using $\|f\|_{s}+\sum_{j=1}^d\|\partial_j f\|_s$ as an equivalent norm in $H^{s+1}$ we get  
\begin{multline*}
\frac{\epsilon^{|\alpha|-1}}{\alpha!}\|[L,\partial^\alpha]u\|_s\leq C_1 \sum_{j=1}^d\sum_{0\not=\gamma\leq\alpha} (C_1\epsilon)^{|\gamma|-1}|\alpha-\gamma+e_j|\epsilon^{|\alpha-\gamma|}\frac{\|\partial^{\alpha-\gamma+e_j}u\|_s}{(\alpha-\gamma+e_j)!}\\
+C_1 \sum_{0\not=\gamma\leq\alpha\atop |\gamma|\leq|\alpha|-1} (C_1\epsilon)^{|\gamma|}\epsilon^{|\alpha-\gamma|-1}\frac{\|\partial^{\alpha-\gamma}u\|_s}{(\alpha-\gamma)!}+(C_1\epsilon)^{|\alpha|}\|u\|_{s}
\end{multline*}
where $e_j$ denotes the $j$-th element of the canonical basis of $\R^d$. Hence we get
\[
\sum_{|\alpha|=N}\frac{\epsilon^{|\alpha|-1}}{\alpha!}\|[L,\partial^\alpha]u\|_s\leq C_2\sum_{1\leq|\tilde{\alpha}|\leq N}
\sum_{|\gamma|=N-|\tilde{\alpha}|+1} (C_2\epsilon)^{N-|\tilde{\alpha}|}|\tilde{\alpha}|\frac{\|\partial^{\tilde{\alpha}} u\|_s}{\tilde{\alpha!}}+(C_2\epsilon)^{N}\|u\|_{s}.
\]
Since the number of multi-indices $\gamma$ of length $N-|\tilde{\alpha}|+1$ is dominated by $2^{N-|\tilde{\alpha}|+d}$ we obtain 
\[
\sum_{|\alpha|=N}\frac{\epsilon^{|\alpha|-1}}{\alpha!}\| [L,\partial^\alpha]u\|_s\,d\sigma\leq C'\sum_{j=1}^N (C'\epsilon)^{N-j} j E^\epsilon_j[u]+(C_2\epsilon)^{N}\|u\|_{s}.
\]
The above constants $C_1,C_2,C'$ depend on the operator $L$ but are independent of $\epsilon, N$. Hence we got
\begin{multline*}
\sum_{|\alpha|=N}\frac{\epsilon(\sigma)^{|\alpha|-1}}{\alpha!}\| [L,\partial^\alpha]u(\sigma)\|_s \leq C'N E^{\epsilon(\sigma)}_N[u(\sigma)] \\ +C'\sum_{j=1}^{N-1} (C'\epsilon(\sigma))^{N-j} j E^{\epsilon(\sigma)}_j[u(\sigma)]+(C_2\epsilon(\sigma))^{N}\|u(\sigma)\|_{s}.
\end{multline*}
Substituting in the left hand side of \eqref{comm} we see that it is sufficient to prove that
\begin{equation}\label{agg0}
\int_0^t \frac{\epsilon(t)^{N-1}}{\epsilon(\sigma)^{N-1}} e^{-C\sigma}(C_2\epsilon(\sigma))^{N}\|u(\sigma)\|_{s}\, d\sigma\leq \frac{1}{8}e^{-Ct}\Phi(t) \qquad \textrm{for} \quad N \geq 1,
\end{equation}
and
\begin{equation}\label{pos2}
\int_0^t \frac{\epsilon(t)^{N-1}}{\epsilon(\sigma)^{N-1}} e^{-C\sigma} C'\sum_{j=1}^{N-1} (C'\epsilon(\sigma))^{N-j} j E^{\epsilon(\sigma)}_j[u(\sigma)]\,d\sigma\leq \frac{1}{8}e^{-Ct}\Phi(t) \qquad \textrm{for} \quad N \geq 2,
\end{equation}
if $\epsilon(0)=\epsilon_0\leq\overline{\epsilon}_0$ with $\overline{\epsilon}_0$ small enough, and the constant $C_0$ in \eqref{tesi} is large enough. \par
Let us prove estimate \eqref{agg0}.  We apply the energy estimate \eqref{energy} with $v=u$. Since $Lu=\mathcal{N}[u]$ has the form in \eqref{ipo10}, \eqref{ipo1bis}, with the aid of the Schauder's estimates $\|fg\|_s\lesssim\|f\|_s\|g\|_{L^\infty}$ one obtains
\[
e^{-C\sigma}\|u(\sigma)\|_s\leq C\|u_0\|_s+C_3\int_0^\sigma (1+\|u(\tau)\|^{p-1}_{L^\infty})e^{-C\tau}\|u(\tau)\|_s\,d\tau
\]
for a suitable constant $C_3>0$. By an application of Lemma \ref{gronwall} with $a(t)\equiv 1$, $\psi(t)=e^{-Ct}\|u(t)\|_s$, $g(t)=C \|u_0\|_s$ and $h(t)= C_3 (1+\|u\|_{L^{\infty}}^{p-1}),$ we get 
\begin{equation} \label{ultimastima}
e^{-C\sigma}\|u(\sigma)\|_s\leq C\|u_0\|_s e^{C_3\int_0^\sigma (1+\|u(\tau)\|^{p-1}_{L^\infty})\,d\tau},
\end{equation} which gives, using $\epsilon(t)^{N-1}/\epsilon(\sigma)^{N-1}\leq 1$ and $\epsilon(\sigma)\leq\epsilon_0e^{-A\sigma}$, 
\begin{align*}
\int_0^t \frac{\epsilon(t)^{N-1}}{\epsilon(\sigma)^{N-1}}& e^{-C\sigma}(C_2\epsilon(\sigma))^{N}\|u(\sigma)\|_{s}\, d\sigma\\
&\leq C\|u_0\|_s e^{C_3\int_0^t (1+\usp)\,d\sigma}\int_0^t (C_2\epsilon_0)^N e^{-AN\sigma}\,d\sigma\\
&\leq \frac{C(C_2\epsilon_0)^N}{AN}\|u_0\|_s e^{C_3\int_0^t (1+\usp)\,d\sigma}\leq \frac{1}{8}e^{-Ct}\Phi(t)
\end{align*}
if the constant $C_0$ in \eqref{tesi} satisfies $C_0\geq C+C_3$ and $C_0\geq C\|u_0\|_s$ (hence these constraints on $C_0$ have to be added to those in \eqref{posizione}), $A\geq 1$ and $\epsilon_0\leq \overline{\epsilon}_0:=1/(8C_2)$. This proves \eqref{agg0}. \par
Let us now prove \eqref{pos2} for $N \geq 2.$ We can use then the inductive hypothesis \eqref{tesi} for subscripts $j<N$ and the estimates $\epsilon(t)^{N-1}/\epsilon(\sigma)^{N-1}\leq e^{-A(t-\sigma)(N-1)}$, as well as $\epsilon(t)\leq \epsilon_0e^{-At}$. 
Then the left-hand side of \eqref{pos2} is
\begin{align*}
&\leq C'C_0 e^{C_0\int_0^t \usp\,d\sigma}\sum_{j=1}^{N-1} (C'\epsilon_0)^{N-j}\int_0^t j e^{(A(j-1)+C_0-C)\sigma-At(N-1)}\,d\sigma\\
&\leq 
e^{-Ct}\Phi(t)  \sum_{j=1}^{N-1} (C'\epsilon_0)^{N-j} \frac{C'j}{A(j-1)+C_0-C}.
\end{align*}
The above sum is $<1/8$ for all $\epsilon_0\leq\overline{\epsilon}_0$ if $\overline{\epsilon}_0$ is small enough, depending on $C$ and $C'$  (unifomely with respect to $C_0>2C$, cf.\ \eqref{posizione}, and $A\geq1$).  \qed

\section{Proof of Proposition \ref{prononlin}}
To simplify the notation and without loss of generality we shall prove Proposition \ref{prononlin} in the scalar case $n=1,$ and for a nonlinear term $\mathcal{N}[u] =g(t,x)u^p$, where $p \in \N, p \geq 2$ and $g(t,x)$ is a smooth function defined on $\overline{\R}_+\times\R$ and satisfying the estimate \eqref{ipo1bis}. 
Let us first consider the case $N=1.$ We have 
\begin{multline*}
E_1^{\epsilon}[\mathcal{N}[u]]= E_1[\mathcal{N}[u]]= \sum_{j=1}^d \|\partial_{x_j}(g(t,x)u^p)\|_s \\ \leq C_1 \sum_{j=1}^d \|\partial_{x_j}g\|_{L^{\infty}(\overline{\R}_+\times\R)} \cdot \|u\|_s^p +C_1 \|g\|_{L^{\infty}(\overline{\R}_+\times\R)} \cdot \|u\|_{L^{\infty}}^{p-1} \cdot E_1[u] \\ \leq C'(\|u\|_s^{p} + \|u\|_{L^{\infty}}^{p-1} \cdot E_1[u]). \end{multline*}
To prove \eqref{nonlin} for $N=1$ we are then reduced to prove that 
$$C'\int_0^t e^{-C\sigma} \|u(\sigma)\|_s^p \, d\sigma \leq \frac{1}{4}e^{-Ct}\Phi(t).$$
Using the estimate \eqref{ultimastima} and assuming  $C_0> p(C_3+C)$ we obtain that
\begin{align*}\int_0^t e^{-C\sigma} \|u(\sigma)\|_s^p\,d\sigma &\leq (C\|u_0\|_s)^p\int_0^t e^{(C(p-1)+C_3p)\sigma} e^{C_3p \int_0^\sigma \|u(\tau)\|_{L^{\infty}}^{p-1}\, d\tau} \\ &\leq (C\|u_0\|_s)^p e^{(C_0-C)t}\int_0^t e^{-(C_0 -Cp-C_3p)\sigma} e^{C_3p \int_0^\sigma \|u(\tau)\|_{L^{\infty}}^{p-1}\, d\tau}\,d\sigma\\ &\leq  \frac{(C\|u_0\|_s)^p}{C_0} e^{-Ct}\Phi(t) \int_0^t  e^{-(C_0 -Cp-C_3 p)\sigma}\, d\sigma \\ &= \frac{(C\|u_0\|_s)^p}{C_0(C_0 -Cp-C_3p)} e^{-Ct}\Phi(t)(1-e^{-(C_0 -Cp-C_3p)t}) \\ &\leq \frac{(C\|u_0\|_s)^p}{C_0(C_0 -Cp-C_3p)} e^{-Ct}\Phi(t). \end{align*} 
Then one can choose $C_0$ large enough to have $\frac{(C\|u_0\|_s)^p}{C_0(C_0 -Cp-pC_3)} \leq \frac{1}{4}.$ The proposition is then proved in the case $N=1.$ \par For $N \geq 2,$ using Leibniz formula we have
$$\| \partial_x^{\alpha}\mathcal{N}[u]\|_s \leq  \|g(t,\cdot)u^{p-1}\partial_x^{\alpha}u\|_s + \hskip-4pt \sum_{\gamma + \alpha_1+\ldots +\alpha_{p}=\alpha\atop |\alpha_j|<|\alpha| \forall j} \frac{\alpha!}{\gamma! \alpha_1! \ldots \alpha_p!} \| \partial_x^{\gamma}g(t,\cdot) \partial_x^{\alpha_1}u \cdot \ldots \cdot \partial_x^{\alpha_p}u \|_s.$$ Then by Schauder's lemma and using the condition \eqref{ipo1bis} we obtain
\begin{multline*} 
E_N^{\varepsilon}[\mathcal{N}[u]] = \sum_{|\alpha|=N}\frac{\epsilon^{|\alpha|-1}}{\alpha!}\| \partial_x^{\alpha}\mathcal{N}[u]\|_s 
 \leq C \| g\|_{L^{\infty}(\overline{\R}_+\times\R)}\cdot  \| u\|_{L^{\infty}}^{p-1} \cdot E_N^{\epsilon}[u]  \\ + C \sum_{|\alpha|=N}\frac{\epsilon^{|\alpha|-1}}{\alpha!} \sum_{\stackrel{\gamma + \alpha_1+\ldots +\alpha_{p}=\alpha}{|\alpha_j|<|\alpha| \forall j}} \frac{\alpha!}{\gamma! \alpha_1! \ldots \alpha_p!} \| \partial_x^{\gamma}g(t,\cdot) \|_{L^{\infty}(\overline{\R}_+\times\R)} \cdot \|\partial_x^{\alpha_1}u\|_s \cdot \ldots \cdot \|\partial_x^{\alpha_p}u\|_s \\ \leq CM \|u\|^{p-1}_{L^{\infty}} \cdot E_N^{\varepsilon}[u] +C' \varepsilon^{p-1}\sum_{|\alpha|=N} \sum_{\stackrel{ \alpha_1+\ldots +\alpha_{p}=\alpha}{|\alpha_j|<|\alpha| \forall j}} (M\epsilon)^{|\gamma|}\prod_{j=1}^p \frac{\varepsilon^{|\alpha_j|-1}\|\partial_x^{\alpha_j}u\|_s}{\alpha_j!} \\ \leq C'  \|u\|^{p-1}_{L^{\infty}} \cdot E_N^{\varepsilon}[u] + C' \varepsilon \left( \sup_{1 \leq j \leq N-1} E_j^{\varepsilon}[u] \right)^p,
\end{multline*}
where we used the fact that $p-1 \geq 1$ and we assume $\epsilon\leq \min\{1,M^{-1}\}$, $M$ being the constant appearing in \eqref{ipo1bis}. 
Then we obtain 
\begin{multline*}
\int_0^t \frac{\epsilon(t)^{N-1}}{\epsilon(\sigma)^{N-1}}e^{-C\sigma}E^{\epsilon(\sigma)}_N[\mathcal{N}[u(\sigma)]]\,d\sigma\\ 
\leq C' \int_0^t \frac{\epsilon(t)^{N-1}}{\epsilon(\sigma)^{N-1}} \usp e^{-C\sigma} E_N^{\epsilon(\sigma)}[u(\sigma)] \, d\sigma \\ + 
C' \int_0^t \frac{\epsilon(t)^{N-1}}{\epsilon(\sigma)^{N-1}} \epsilon (\sigma) e^{-C\sigma}  \left( \sup_{1 \leq j \leq N-1} E_j^{\epsilon(\sigma)}[u(\sigma)] \right)^p \, d\sigma
\end{multline*}
if $\epsilon\leq \overline{\epsilon}_0:=\min\{1,M^{-1}\}$.\par
To conclude the proof, it is then sufficient to prove that 
\begin{equation} \label{pluto}
 \int_0^t \frac{\epsilon(t)^{N-1}}{\epsilon(\sigma)^{N-1}} \epsilon (\sigma) e^{-C\sigma}  \left( \sup_{1 \leq j \leq N-1} E_j^{\epsilon(\sigma)}[u(\sigma)] \right)^p \, d\sigma \leq \frac{1}{4}e^{-Ct} \Phi(t).
\end{equation}
First of all we observe that since the function $\epsilon (t)$ is decreasing, we have $ \frac{\epsilon(t)^{N-1}}{\epsilon(\sigma)^{N-1}} \leq 1$ for $0\leq \sigma \leq t.$ Moreover, writing
$$e^{-C\sigma} \left( \sup_{1 \leq j \leq N-1} E_j^{\epsilon(\sigma)}[u(\sigma)] \right)^p = e^{(p-1)C\sigma} \left( \sup_{1 \leq j \leq N-1} e^{-C\sigma}E_j^{\epsilon(\sigma)}[u(\sigma)] \right)^p$$
and using the inductive assumption \eqref{tesi}, we obtain that
$$ \left( \sup_{1 \leq j \leq N-1} E_j^{\epsilon(\sigma)}[u(\sigma)] \right)^p \leq C_0^p e^{Cp \int_0^t \usp \, d\sigma}. $$
Hence, using the explicit expression of $\epsilon(t)$ in \eqref{epst} and choosing $A >C_0(p-1) \geq C(p-1)$, we get
\begin{multline*}
\int_0^t \frac{\epsilon(t)^{N-1}}{\epsilon(\sigma)^{N-1}} \epsilon (\sigma) e^{-C\sigma}  \left( \sup_{1 \leq j \leq N-1} E_j^{\epsilon(\sigma)}[u(\sigma)] \right)^p \, d\sigma \\ \leq C_0^p \epsilon_0 \int_0^t e^{((p-1)C-A)\sigma  -(A-C_0p)\int_0^{\sigma}\|u(\tau)\|^{p-1}_{L^{\infty}}\, d\tau }\, d\sigma \\ \leq C_0^p \epsilon_0 e^{C_0 \int_0^t \usp \, d\sigma} \cdot \int_0^t e^{-(A-(p-1)C)\sigma \,d\sigma}  \\ = \frac{C_0^p \epsilon_0}{A-(p-1)C}(1-e^{-(A-(p-1)C)t}) e^{C_0 \int_0^t \usp \, d\sigma} \\ \leq \frac{C_0^{p} \epsilon_0}{A-(p-1)C}e^{-Ct} \Phi(t) \leq \frac{1}{4}e^{-Ct}\Phi(t)
\end{multline*}
choosing $A$ sufficiently large. This concludes the proof. \qed

\section{Examples and remarks}\label{osservazioni}
The following two examples make clear the expression for $\epsilon(t)$ given in \eqref{epst}.
\begin{example}\rm 
 Consider first the linear Cauchy problem (hence $p=1$) in $\overline{\R}_+\times\R$ given by
\[
\begin{cases}
\partial_tu-x\partial_x u=0\\
u(0,x)=u_0(x)=(1+x^2)^{-1}.
\end{cases}
\]
The initial datum extends to a holomorphic function $u_0(x+iy)$ in the strip $|y|<1$, whereas the solution, given by $u(t,x)=(1+x^2e^{2t})^{-1}$, extends, for every fixed $t$, to a holomorphic function in the strip $|y|<e^{-t}$. Hence the linear part of the equation is already responsible of the exponential decay in \eqref{epst}. For equations with constant coefficients  (even non-linear) one instead expects algebraic decay for the radius of analyticity (cf. e.g.\ \cite{bo2}). We plan to present a systematic study of this special case in a future work.\par
\end{example}

\begin{example} \rm
Fix now $p\geq 2$, $p\in\N$, and consider instead the non-linear Cauchy problem
\[
\begin{cases}
\partial_tu=\frac{1}{p-1}u^p\\
u(0,x)=u_0(x)=(1+x^2)^{-1}.
\end{cases}
\]
The maximal solution is defined in $[0,1)\times\R$ and reads 
\[
u(t,x)=((1+x^2)^{p-1}-t)^{-1/(p-1)}.
\] The initial datum is the same as before. Now the solution extends to a holomorphic function in the strip $|y|<\sqrt{1-t^{1/(p-1)}}\sim(p-1)^{-1/2}\sqrt{1-t}$ as $t\to 1$. This agrees with the formula \eqref{epst} with $A=1/2$; in fact, for $0\leq t<1$,
\[
\exp\big(-\frac{1}{2}\int_0^t \|u(\sigma)\|^{p-1}_{L^\infty}\,d\sigma)\big)=\exp\big(-\frac{1}{2}\int_0^t (1-\sigma)^{-1}\,d\sigma\big)=\sqrt{1-t}. 
\]
\end{example}

\begin{remark}\rm  A careful inspection of the proof shows that the conclusions of Theorems \ref{mainteo} and \ref{mainteo2} hold with a radius of analyticity  
\[
\epsilon(t)=\frac{1}{B}\exp\big(-A\int_0^t (1+\|u(\sigma)\|^{p-1}_{L^\infty})\,d\sigma \big)
\]
for a constant $A>0$ depending only on the equation, and a constant $B$ depending on the initial datum; to be precise $B\simeq (\|u_0\|_{s}+\sup_{N\geq 1}E^{\epsilon_0}_{N}[u_0])^{p-1}$. \par 
\end{remark}
\begin{remark} \label{strictly}\rm
We observe that the condition \eqref{simbolo2} is used in the proof of Theorem \ref{mainteo2} only to prove the energy estimate \eqref{energy}. Hence Theorem \ref{mainteo2} still holds for every system admitting a smooth solution satisfying an estimate of the form \eqref{energy}. In particular, we can apply our result to a strictly hyperbolic system of the form $$L'u =\mathcal{N}[u]$$ $\mathcal{N}[u]$ is of the form \eqref{ipo10}, \eqref{ipo1bis} and 
$$L'= \partial_t + A(t,x,D_x)$$ for some $n\times n$ matrix $A(t,x,\xi)$ satisfying the condition \eqref{simbolo} and admitting distinct purely imaginary eigenvalues $i\lambda_j(t,x,\xi), j=1,\ldots, n$, with
$$\lambda_1 (t,x,\xi)<\lambda_2 (t,x,\xi)<\ldots <\lambda_n (t,x,\xi), \quad (t,x,\xi) \in \overline{\R}_+ \times \rd \times \rd.$$
It is well known that under these assumptions for $s$ sufficiently large the Cauchy problem $$\begin{cases}L'u=\mathcal{N}[u] \\ u(0,x)=u_0 \end{cases}$$ admits a solution $u \in C^1(I,H^s)\cap C^0(I,H^{s+1})$ defined for $(t,x) \in I \times \rd$ for some time interval $I$ and the energy estimate \eqref{energy} holds true. Hence the arguments in the proof of Theorem \ref{mainteo2} can be applied to prove the analyticity of the solution and the estimate of the radius of analyticity.
\end{remark}

\end{document}